\documentclass[10pt]{article}
\usepackage{amsmath, exscale, epsfig,amssymb, 
 makeidx}

\small
\normalsize
\usepackage[latin1]{inputenc}
\usepackage[T1]{fontenc}
\usepackage{amsfonts}
\usepackage{amsbsy}
\usepackage{amscd}
\usepackage{theorem}

\usepackage{epsf}
\usepackage{psfrag}
\usepackage{graphicx}

\usepackage{multicol}
\allowdisplaybreaks

\def\build#1_#2^#3{\mathrel{\mathop{\kern 0pt#1}\limits_{#2}^{#3}}}
\def\noi{{\noindent}}

\def\cq{$\hfill \square$}
\def\un{{\bf 1}}

\newcommand{\bE}{{\bf E}}

\newcommand{\bU}{{\bf U}}
\newcommand{\buu}{{\bf u}}
\newcommand{\bvv}{{\bf v}}

\newcommand{\bN}{\mathbb{N}}

\newcommand{\bP}{{\bf P}}

\newcommand{\bbU}{\mathbb{U}}

\newcommand{\cB}{\mathcal{B}}

\newcommand{\cF}{\mathcal{F}}
\newcommand{\cG}{\mathcal{G}}
\newcommand{\cH}{\mathcal{H}}

\newcommand{\cM}{\mathcal{M}}

\newcommand{\cP}{\mathcal{P}}

\newcommand{\cT}{\mathcal{T}}

\def\ovvm{{\rm m}}

\def\be{\begin{equation}}
\def\boU{ \partial \bbU}
\def\ee{\end{equation}}
\def\ba{\begin{eqnarray*}}
\def\ea{\end{eqnarray*}}

\def\noi{\noindent}

\def\supp{{\rm supp\,}}

\def\cqfd{ \hfill $\blacksquare$ }
\newcommand{\bbT}{\mathbb{T}}
\newcommand{\cut}{{\rm Cut}}

\newtheorem{theorem}{Theorem}[section]
\newtheorem{lemma}[theorem]{Lemma}
\newtheorem{proposition}[theorem]{Proposition}

\newtheorem{definition}{Definition}[section]
{\theorembodyfont{\rmfamily}}
{\theorembodyfont{\rmfamily}\newtheorem{remark}{Remark}[section]}
\begin{document}

\title{ {\bf AN ELEMENTARY PROOF OF HAWKES'S CONJECTURE  
ON GALTON-WATSON TREES.}}
\author{Thomas {\sc Duquesne}  
\thanks{Laboratoire de Probabilit\'es et Mod\`eles Al\'eatoires; 
Universit\'e Paris 6, 16 rue Clisson, 75013 PARIS, FRANCE.  Email: thomas.duquesne@upmc.fr }}

\vspace{2mm}
\date{\today} 

\maketitle

\begin{abstract}
In 1981, J. Hawkes conjectured the exact form of the Hausdorff gauge function for the boundary of supercritical Galton-Watson trees under a certain assumption on the tail at the infinity of the total mass of the branching measure. Hawkes's conjecture has been proved by T. Watanabe in 2007 as well as other other precise results on fractal properties of the boundary of Galton-Watson trees. The goal of this paper is to provide an elementary proof of Hawkes's conjecture under a less restrictive assumption than in T. Watanabe's paper, by use of size-biased Galton-Watson trees introduced 
by Lyons, Pemantle and Peres in 1995. \\

\end{abstract}

\section{Introduction.}
\label{introd} Fractal properties of the boundary of supercritical Galton-Watson trees have been intensively studied since the seminal paper by R.A. Holmes \cite{Holmes73}, who first studied the exact Hausdorff measure of a specific case, and since the paper by J. Hawkes \cite{Hawkes81}  who determined the growth number of a Galton-Watson tree in the general case and proved 
that the Hausdorff  dimension of its boundary is the logarithm of the mean of the offspring distribution (see also Lyons \cite{Ly90} for a simple proof). The problem of finding an exact Hausdorff function in the general setting has been studied by Q. Liu  \cite{Liu96} who considered a large class of offspring distributions. Packing measure, thin and thick points as well as multifractal properties of the branching measure have been also investigated by Q. Liu \cite{Liu01, Liu00}, P. Mörters and N. R. Shieh \cite{MorShieh02, MorShieh04} and T. Watanabe \cite{TWata04}. 

 In the cases studied by Q. Liu \cite{Liu96}, the corresponding Hausdorff measure coincides with the branching measure. This was predicted by J. Hawkes \cite{Hawkes81} who conjectured the general form of the Hausdorff gauge function for the boundary of supercritical Galton-Watson trees under a natural assumption on the right tail of the distribution of the total mass of the branching measure. This long standing conjecture has been recently solved by  T. Watanabe in \cite{TWata07}: in this paper, T. Watanabe provides necessary and sufficient conditions for the 
existence of an exact Hausdorff measure that is absolutely continuous with respect to the branching measures as well as precise results on several important examples. The goal of our paper is to provide an elementary proof of Hawkes's conjecture that holds true under a less restrictive assumption than in Watanabe's paper and that relies on size-biased 
Galton-Watson trees that have been introduced in \cite{LyoPemPer} by R. Lyons, R. Pemantle and Y. Peres.

  Let us briefly state Hawkes's conjecture (we refer to Section \ref{basics} for more formal definitions).  To simplify notation, we assume that all the random variables that we consider are defined on the same probability space $(\Omega, \cF , \bP)$ that is assumed to be complete and sufficiently large to carry as many independent random variables as we need. Let 
$\xi= (\xi (k), k \in \bN)$ be a probability distribution on $\bN$ that is viewed as the offspring distribution of a Galton-Watson tree $\cT $. Informally, $\cT$ is the family-tree of a population stemming from one ancestor and evolving randomly as follows: each individual of the population has an independent random number of children distributed in accordance with the offspring distribution $\xi$. We assume that the ancestor is at generation $0$, its children are at generation $1$, their children are at generation $2$... and so on; if we denote by $Z_n (\cT)$ the number of individuals at generation $n$, then $(Z_n (\cT)\, ; \, n \geq 0)$ is a Galton-Watson Markov chain with offspring distribution $\xi$ that starts at state $1$. Let us set 
\begin{equation}
\label{notationm}
 \ovvm= \sum_{k \in \bN} k \, \xi (k) \in [0, \infty] \quad {\rm and} \quad f (r) =  \sum_{k \in \bN}   \xi (k) r^k \; , \; r \in [0, 1] \; .
 \end{equation}
The function $f$ is the generating function of $\xi$. It is convex and the equation $f(r)=r$ has at 
most two roots in $[0, 1]$. We denote by $q$ the smallest one ($1$ being obviously the largest one). Standard results on Galton-Watson Markov chains imply that 
$\bP ( \# \cT < \infty )= q $. Morever $q = 1$ iff $\ovvm \leq 1$. Therefore, if $\ovvm >1$, $\bP (\# \cT = \infty) >0$.  
We shall make the following assumptions on $\xi$: 
\begin{equation*}
\label{assumption}
\tag{ Hyp($\xi$) }  \ovvm \in (1, \infty)  \qquad {\rm and}\qquad \sum_{k \geq 2} k \log (k)\,  \xi (k) \, < \, \infty \; .
\end{equation*}
We suppose $\ovvm >1$ because we want to study $\cT$ on the event $\{ \# \cT = \infty \}$. The second assumption on $\xi$ is motivated by the following result known as Kesten-Stigum's Theorem that asserts that, under Hyp($\xi$), we have 
\begin{equation}
\label{KSconv}
 \lim_{n \rightarrow \infty} \frac{Z_n (\cT) }{\ovvm^n} = W \quad  {\rm a.s.} \; \,  {\rm and} \;\,   {\rm in} \; {\rm L}^1 (\Omega , \cF , \bP) \; .
\end{equation}
Consequently, $W$ has unit expectation: $\bE [W]=1$. Moreover,  $\un_{\{ W>0 \}}= \un_{ \{ \# \cT= \infty  \} } $ almost surely and the exact rate of growth of $\cT$ is $n \mapsto \ovvm^n W$. 
We refer to the original paper by Kesten and Stigum \cite{KesStig} and to the paper by R. Lyons, R. Pemantle and Y.  Peres  \cite{LyoPemPer} for an elementary proof of (\ref{KSconv}).

  We view $\cT$ as an ordered rooted tree. Namely, we take the ancestor of the population as the root and we associate a rank of birth with any individual. In this way we can label each individual by a finite word 
of positive integer $(i_1, \ldots, i_n)$: the length of the word $n$ is the generation of the labelled individual and $i_k$ represents the birth-rank of its ancestor at generation $k$. 
The set of words created in this way completely encodes $\cT$ which can be therefore viewed  
as a random subset of the set of finite words written with positive integers. We denote by $\bbU$ the set of finite words written by positive integers. We are interested in the infinite boundary of $\cT$ that is the set of all infinite lines of descent in $\cT$. We denote by $\partial \cT$ the boundary of $\cT$. Let us assume that $\partial \cT$ is non-empty (which is equivalent to assume that $\# \cT $ is infinite); we denote by $\bN^*= \bN \backslash \{ 0 \}$ the set of positive integers; to each infinite line of descent we associate an infinite $\bN^*$-valued sequence $(i_n\, ; \, n \geq 1)$ such that for any $n \geq 1$, the finite word  $(i_k\, ; \, 1\leq k \leq n)$ is the label of the individual at generation $n$ belonging to the infinite line of descent. Therefore, we can view $\partial \cT$ as a random subset of the set of $\bN^*$-valued and $\bN^*$-indexed sequences that we denote by $\boU$. We equip $\boU$ with the 
metric $\delta$ defined as follows: if $\buu = (i_k\, ; \, k \geq 1)$ and 
$\bvv = (j_k\, ; \, k \geq 1)$ are two elements of $\boU$, then $\delta (\buu, \bvv)= \exp (-n)$, 
where $n$ is the largest integer $m$ such that $ \buu$ and $\bvv$ agree on the first $m$ terms  
($n$ is taken as $0$ if $i_1\neq j_1$). The resulting metric space $(\boU, \delta)$ is complete and separable and $\partial \cT$ is actually a random compact subset of $\boU$. 

   Hawkes's conjecture concerns the problem of finding a gauge function $g$ such that the $g$-Hausdorff measure of $\partial \cT$ is positive and finite. Before stating the conjecture, 
let us briefly recall standard definitions about Hausdorff measures: we restrict our attention to $g$-Hausdorff measures with sufficiently regular gauge function $g$; more precisely, we say that $g$ is a regular function if firstly, there exists an interval $(0, r_0)$ on which $g$ is right-continuous non-decreasing, secondly, $\lim_0 g = 0$ and thirdly, there exists $C >1$ such that 
$g(2r) \leq C g(r)$, for any $r \in (0, r_0/2)$ (this last assumption is often called the "doubling condition" though some authors use the term of  "blanketed" Hausdorff function). Then, the $g$-Hausdorff measure on $(\boU, \delta)$ is defined as follows: for any $A \subset \boU$ and for any $\varepsilon \in (0, r_0) $, we first set 
$$  \cH^{(\varepsilon)}_g \left( A \right) =  \inf \left\{  
\sum_{ n\in \bN} g \left( {\rm diam} (C_n) \right) \; \,  ; \; 
A \subset \bigcup_{n \in \bN} C_n  \; \; {\rm and} \; \; {\rm diam} (C_n)
 \leq \varepsilon \, , \, n \in \bN  \right\} , $$
  where ${\rm diam} (C)= \sup \{ \delta (\buu, \bvv) \, ; \, \buu, \bvv \in C \}$ stands for the diameter of a subset $C \subset \boU$; then the $g$-Hausdorff measure is given by $ \cH_g \left( A \right) = \lim_{\varepsilon \rightarrow 0} \; \,  \cH^{(\varepsilon)}_g \left( A \right) \; \in [0, \infty ]$.

\medskip

\noi $\bullet$ {\it Hawkes's conjecture.} Let $\xi$ be a probability measure on $\bN$ that satisfies Hyp($\xi$). Let $\cT$ be a Galton-Watson tree with offpring distribution $\xi$. Let $W$ be defined by (\ref{KSconv}). We set 
\begin{equation}
\label{defFF}
 F(x):= -\log \bP ( W > x ) \; .
\end{equation} 
We first assume that $F$ is regularly varying at $\infty$; more precisely, we suppose that $F$ is of the following form: 
\begin{equation}
\label{Hawkes}
 F(x)= x^b \ell (x)
\end{equation} 
where $b >0$ and $\ell$ is slowly varying function at $\infty$. We denote by 
$F^{-1}$ the right-continuous inverse of $F$: $ F^{-1} (x) = \inf \{ y \geq 0 \; : \; F(y) > x\,      \}  $
and we set 
\begin{equation}
\label{defgg}
g(r) = r^{\log \ovvm} \, F^{-1}\!  \left(\,  \log \log 1/r  \, \right) \; , \quad r \in (0, e^{-1} ) \; .
\end{equation}
Then, Hawkes \cite{Hawkes81} p.382 conjectured that under (\ref{Hawkes}), there exists $c_{\xi} \in (0, \infty)$ that only relies on $\xi$ such that 
\begin{equation}
\label{trueconj}
 \cH_g \left( \partial \cT \right) = c_{\xi} W \; .
 \end{equation}
As already mentioned, this result has been solved by T. Watanabe (2007) in Theorem 1.6 \cite{TWata07} under the following assumption that is weaker than (\ref{Hawkes}): for any sufficiently large $x$ 
\begin{equation}
\label{Wata}
 A^{-1}    x^b \ell (x) \leq F(x) \leq A  x^b \ell (x)  \; , 
\end{equation}
where $A$ is a constant larger than one. The paper by T. Watanabe \cite{TWata07} contains other general results. (Actually, Watanabe's proof of Hawkes's conjecture is a consequence of Theorem 1.2 \cite{TWata07} that provides a general criteria  to decide whether the branching measure is an exact Hausdorff measures with a regular gauge function.)

  The goal of this paper is to give an alternative short proof of  Hawkes's conjecture under a less restrictive assumption than 
(\ref{Wata}) and by use of different techniques that we claim to be elementary. Before stating the main result of this paper, let us mention that the branching measure $M$, whose definition is recalled in Section \ref{basics}, is a finite random measure on $\boU$ 
associated with $\partial \cT$  and such that $M (\partial \cT)= W$.
\begin{theorem}
\label{mainres}
Let $\xi$ be a probability measure on $\bN$ which satisfies Hyp($\xi$). Let $\cT$ be a Galton-Watson tree with offpring distribution $\xi$. Let $W$ be defined by (\ref{KSconv}); let $F$ be defined by (\ref{defFF}) and let $g$ be defined by (\ref{defgg}). We assume 
\begin{equation}
\label{weakHawkes}
\sup_{ x \in [1, \infty )} \frac{F^{-1} (2x)}{F^{-1} (x)} < \infty \; .
 \end{equation} 
Then, there exists $c_{\xi} \in (0, \infty) $ that only depends on $\xi$ 
such that 
\begin{equation}
\label{typicaluppdens}
 \bP\! {\rm -a.s.} \; \,  {\rm for} \; M \! {\rm -almost} \; {\rm all} \; \buu \; : \quad \limsup_{r \rightarrow 0} \frac{M (B(\buu , r)) }{g(r)} = c_{\xi}^{-1} \; , 
\end{equation} 
where $B(\buu , r)$ stands for the open ball in $(\boU, \delta)$ with center $\buu$ and radius $r$. Furthermore, we have 
 \begin{equation}
\label{exactH}
\bP{\rm -a.s} \qquad \cH_g \left( \, \cdot \, \cap \partial \cT \right) = c_{\xi} \cdot  M \; ,
\end{equation} 
where $M$ stands for the branching measure associated with $\partial \cT$. 
\end{theorem}
  
\begin{remark}
\label{mildreg}
Assumption (\ref{weakHawkes}) is weaker than (\ref{Hawkes}). More precisely, it is easy to prove that  (\ref{weakHawkes}) is equivalent to the following: 
\begin{equation}
\label{recipwH} 
 \exists a >0 \; , \quad F^{-1}(sx) \leq 2^a  s^{a} F^{-1}(x) \; , \quad s , x \in [1, \infty)   \; .
\end{equation}
Therefore, $F$ satisfies 
\begin{equation}
\label{conseqF} 
 2^{-1} s^{1/a}F(x) \leq  F (sx) \; , \quad s \geq 2^a  , x \geq F^{-1}(1)   \; .
\end{equation}
$\; $  \cq 
\end{remark}
\begin{remark}
\label{hyptranslation}
There is no known necessary and sufficient condition expressed in terms of $\xi$ for $F$ to satisfy (\ref{weakHawkes}) (neither for (\ref{Wata}) nor (\ref{Hawkes})). However specific cases have been considered by Q. Liu and T. Watanabe: see \cite{Liu96} and \cite{TWata07}. Let us also mention that when (\ref{exactH}) holds true, there is no simple general closed formula giving $c_\xi$ in terms of $\xi$. The avaible results characterizing $c_\xi$ are either quite involved or they require the knowledge of the distribution of $W$: see Theorem 1.1 in \cite{TWata07} for a general characterization of $c_\xi$;  see Liu (Theorem 1\cite{Liu96}) or Watanabe (Theorem 1.6  \cite{TWata07})when the support of $\xi$ is bounded; see also Watanabe (Theorem 1.6 \cite{TWata07})  when $F$ satisfies (\ref{Wata}). \cq 
\end{remark}

\section{Notation and basic definitions.}
\label{basics}
Let us start with basic notation: we denote by $\bN$ the set of nonnegative integers and by $\bN^*$ the set of positive integers; let $\bbU= \bigcup_{n \in \bN} (\bN^*)^n $ be the set of finite words written with positive integers, with the convention $  (\bN^*)^0= \{ \varnothing\}$. Let $u=(i_k \, ; \,  1\leq k \leq n) \in \bbU$; we set 
$|u|=n$ that is the length of $u$, with the convention $| \varnothing  |= 0$. Words of unit length are identified with positive integers. For any $m \in \bN$ we set 
$u_{|m}= (i_k \, ; \, 1 \leq k \leq n\wedge m)$, with the convention $u_{|0}=  \varnothing $; observe that $u_{|m}= u$ if $m \geq n$. Let 
$v=(j_k ;  1\leq k \leq m) \in \bbU$, we define $u*v \in \bbU$ by the word $(\ell_k \,  ;  \, 1\leq k \leq n+m)$ where $\ell_k= i_k$ if $k \leq n$ and $\ell_k = j_{k-n}$ if $ k >n $: the word $u*v$ is the concatenation of $u$ and $v$ (observe that $\varnothing * u= u * \varnothing=u $). We next introduce the genealogical order $\preceq $ by writing 
$u \preceq v$ iff  $v_{|_{ |u|}}= u$. For any $u,v \in \bbU$, we denote by $u\wedge v$ the $\preceq$-maximal word $w$ such that $w \preceq u$ and $w \preceq v$. For any $u \in \bbU$, we denote by $\bbU_u$ the $\preceq $-successors of $u$. Namely, $\bbU_u$ is the set of words 
$u*v$ where $v$ varies in $\bbU$. On $\bbU_u$, we define the $u$-shift $\theta_u$ by 
$\theta_u (u*v) =v$. 
\begin{definition}
\label{treedef}
A subset $T \subset \bbU$ is a tree iff it satisfies the following conditions. 
\begin{itemize}
\item{{\rm Tree(1)}:}  If $u \in T $, then $ u_{|m} \in T$, for any $m \in \bN$ (in particular $ \varnothing $ belongs to $T$). 
\item{{\rm Tree(2)}:}  For any $u\in \bbU$, there exists $k_u (T)\in \bN \cup \{-1 \}$ such that  the following properties hold true. 
\begin{itemize}
\item If $k_u (T)=-1$, then $u \notin T$.
\item If $k_u (T) =0$, then $ T \cap \bbU_u = \{ u \}$.  
\item If $k_u (T) \geq 1$, then the set of words $ \left\{   u*i \, ; \,   1 \leq i \leq  k_u (T)  \right\}$ is exactly the set of words $v \in T$ such that $ |v|= |u| +1$ and $u \preceq v $.  \cq 
\end{itemize}
\end{itemize}
\end{definition}   
We denote by $\bbT$ the class of subsets of $\bbU$ satisfying Tree(1) and Tree(2). More precisely, this definition provides a canonical coding of finite-degree ordered rooted trees. For sake of simplicity, any element $T$ in  $\bbT$ shall be called a tree.

  For any $m \in \bN$ and any $T \in \bbT$, we set $ T_{|m}= \{ u \in T \; : \; |u| \leq m \} $. Observe that $T_{|m}$ is a finite tree. For any word $u \in \bbU$ and any tree $T$, we define the $u$-shift of $T$ by 
$$ \theta_u T= \theta_u \left( T \cap \bbU_u \right)= \{ v \in \bbU \; : \; u*v \in T \, \} .$$
We see that $\theta_u T$ is empty iff $u \notin T$; in any case, $\theta_u T$ is a tree. 
For any $u \in T$, we define the tree $T$ cut at vertex $u$ as the following subset of $\bbU$: 
$$ \cut_u T = T \backslash \left\{ u* v \; ; \; v \in  \theta_u T \backslash \{ \varnothing \} \, \right\} . $$  
Observe that $u \in \cut_u T$ and that $\cut_u T$ is a tree. Next, for any $T \in \bbT$ and any $n \in \bN$, we set 
$$ Z_n (T) = \# \{ u \in T \; : \; |u|= n \, \} \;  \in \bN .$$
In the graph-terminology $Z_n (T)$ is the number of vertices of $T$ at distance $n$ from the root. If we view $T$ as the family tree of a population whose $ \varnothing $ is the ancestor and whose genealogical order is $\preceq$, then $Z_n (T)$ is the number of individuals at the $n$-th generation.

\medskip 

  We denote by $\boU$ the set $(\bN^*)^{\bN^*}$ of the $\bN^*$-valued and $\bN^*$-indexed sequences. Let $\buu = ( i_k \, ; \, k \geq 1)$ be in $\boU$. For any $m \geq 0$, we set 
 $\buu_{|m}= (i_k\, ; \, 1\leq k \leq m) \in \bbU$, with convention $\buu_{|0}= \varnothing$. If $\bvv \in \boU$, then we denote by $\buu \wedge \bvv$ the $\preceq$-maximal finite word $w$ such that $\buu_{ |_{|w|}}=\bvv_{|_{|w|}}= w$. 
We equip $\boU$ with the following ultrametric $\delta$ given by 
$$ \delta (\buu , \bvv)= \exp \left( -| \buu \wedge \bvv| \right) \; .$$
The resulting metric space $( \boU , \delta)$ is separable and complete and we denote by $\cB (\boU)$ its Borel sigma-field. For any $r \in (0, \infty)$ and for any $\buu \in \boU$, we denote by $B (\buu , r)$ the open $\delta$-ball with center $\buu$ and radius $r$. 
We shall often use the notation 
\begin{equation}
\label{notationnr}
n(r) = \lfloor (-\log (r) \, )_+ \rfloor  +1\; , 
\end{equation}
where $(\cdot )_+$ stands for the positive part function and $\lfloor \cdot \rfloor $ for the integer part function. Observe that $B (\buu , r)= \boU$ if $r >1$; if $r \in (0, 1]$,  
$B (\buu , r)$ is the set of $\bvv$ such that $\bvv_{|n (r) }= \buu_{|n(r)}$. This has several consequences. {\it Firstly}, any open ball is also a closed ball. {\it Secondly}, we have 
$B (\bvv , r)= B (\buu , r)$ for any $\bvv $ in $B (\buu , r)$; therefore there is only a countable number of balls with positive radius; more precisely for any $u \in \bbU$, we set 
\begin{equation}
\label{Budef}
B_u = \left\{ \bvv \in \boU \; : \; \bvv_{|_{|u|}}= u\, \right\} . 
\end{equation}
Then, 
$$ \left\{ B( \buu , r )  \; ; \; r \in (0, \infty) \, , \, \buu \in \boU \,  \right\} = \left\{ B_u \; ;  \; u \in \bbU  \right\}  \; .
$$
{\it Thirdly}, for any pair of balls either there are disjoint or one is contained in the other. We shall further refer to these properties as to the {\it specific properties of balls in $\boU$.}

    Let $T \subset \bbT$. We define the boundary $\partial T$ by the set 
$ \left\{\,  \buu \in \boU \; : \; \buu_{|n} \in T \, , \, n \in \bN^*  \right\} $
Obviously, $\partial T$ is empty iff $T$ if finite; moreover, since $T_{|m}$ is finite for any $m \in \bN$, an easy diagonal-extraction argument implies that $\partial T$ is a compact set of $(\boU, \delta)$.

\medskip
   
   We equip $\bbT$ with the sigma field  $\cG$ generated by the subsets 
\begin{equation}
\label{defineAu}
A_u:= \left\{ T \in \bbT \; : \; u \in T  \right\} \; , \quad u \in \bbU \; .
\end{equation}
 For any $u \in \bbU$, it is easy to check that the application $T \mapsto k_u (T)$ is $\cG$-measurable. Moreover $T \mapsto \theta_u T$ is $(\cG, \cG)$-measurable and for any 
$n \in \bN$, $T \mapsto Z_n (T)$ is $\cG$-measurable.

    Recall that we suppose throughout the paper that all the random variables we need are defined on the same probability space $(\Omega, \cF , \bP)$. Then, a random tree $\cT$ is an 
application from $\Omega$ to $\bbT$ that is $(\cF, \cG)$-measurable. We shall concentrate our attention on a particular class of random trees called Galton-Watson trees that can be recursively defined as follows. 
\begin{definition}
\label{GWdef}
Let $\xi= (\xi (k)\, ; \, k \in \bN)$ be a probability on $\bN$. A random tree $\cT$ is said to be a Galton-Watson tree with offspring distribution $\xi$ (a GW($\xi$)-tree for short) if its distribution on $(\bbT, \cG)$ is characterized by the following conditions. 
\begin{itemize}
\item{{\rm GW(1)}:} the $\bN$-valued random variable $k_{ \varnothing } (\cT)$ is distributed in accordance with $\xi$. 
\item{{\rm GW(2)}:} If $\xi (k) >0$, then under $\bP \left( \, \cdot \, | k_{ \varnothing } (\cT) = k \right)$, 
the first generation subtrees $\theta_1 \cT, \ldots , \theta_k \cT $ are independent with the same distribution as $\cT$ under $\bP$. \cq 
\end{itemize}
\end{definition}
 
\noi This recursive definition induces a unique distribution on $(\bbT, \cG)$:  we refer to Neveu \cite{Ne} for a proof. 

\medskip

  Recall that we assume Hyp($\xi$) and recall (\ref{KSconv}). Then, for any $u \in \bbU$ we set $ W_u:= \limsup \ovvm^{-n} Z_n (\theta_u \cT) $. Definition \ref{GWdef} easily implies that for any $u\in \bbU$, the random variables $\left( W_{u*v} \, ; \, v \in \bbU \right) $ under 
$\bP ( \, \cdot \,  |\,  u \in \cT)$ have the same distribution as the random variables $\left( W_{v} \, ; \, v \in \bbU \right)$ under $\bP$. 
Therefore, we almost surely have  
\begin{equation}
\label{pourtoutu} 
\forall \, u \in \bbU \, \quad W_u= \lim_{n \rightarrow \infty} \, \ovvm^{-n} Z_n (\theta_u \cT) < \infty \quad {\rm and} \quad \un_{\{ W_u>0 \}}= \un_{ \{ 
\partial (\theta_u  \cT )\neq \emptyset \}}. 
\end{equation}
Moreover, since for any $n\geq 1$ and for any $u \in \bbU$, 
$Z_{n} (\theta_u \cT)$ is the sum of the $Z_{n-1} (\theta_{u*i} \cT)$'s over $i \in \bN^*$, we a.s. have 
\begin{equation}
\label{projective}
 \forall u \in \bbU \; , \quad W_{u} = \ovvm^{-1}  \sum_{i \in \bN^*}W_{u*i} \; .  
 \end{equation}
Denote by $\cM_f ( \boU)$ the set of finite measures on 
$( \boU, \cB (\boU))$; equip $\cM_f ( \boU)$ with the topology of weak convergence. Then, a random finite measure on $( \boU, \cB (\boU))$ is an application from $\Omega$ to $\cM_f ( \boU)$ that is measurable with respect to $\cF$ and to the Borel sigma-field  of $\cM_f ( \boU)$. Thanks to (\ref{projective}), the collection of random variables $(W_u\, ; \, u \in \bbU)$ allows to define a random finite measure $M$ on $( \boU, \cB (\boU))$  that is characterized by the following properties. 
\begin{itemize}
\item{{\rm BM(1)}:} Almost surely, $M$ is diffuse and its topological support is $\partial \cT$.  
\item{{\rm BM(2)}:} Almost surely, for any $u \in \bbU$, $M (B_u)= \ovvm^{-|u|} W_u $, where we recall notation $B_{u}$ from (\ref{Budef}).  
\end{itemize}
This define the branching measure associated with $\cT$. It is in some sense the most spread out measure on $\partial \cT$; it is a natural candidate to be a Hausdorff measure on $\partial T$. Actually, Theorem \ref{mainres} is proved by applying to $M$ the following Hausdorff-type comparison results of  measures. 
\begin{proposition}
\label{comparison}
Let $\mu \in \cM_f (\boU)$. Let $g: (0, r_0)\rightarrow (0, \infty)$ be a right-continuous and non-decreasing application such that $\lim_0 g =0$ and such that there exists $C>1$ that satisfies $g(2r)\leq C g(r)$, for any $r\in (0, r_0/2)$. Then, for any Borel subset $A\subset \boU$, the following assertions hold true. 
\begin{itemize} 
\item{(i)} If  $ \limsup_{r \rightarrow 0} \frac{\mu (B(\buu, r))}{g(r)} \leq 1 $ for any $\buu \in A$, then 
$  \cH_g (A)  \geq  C^{-1} \mu (A) $. 
\item{(ii)}  If  $ \limsup_{r \rightarrow 0} \frac{\mu (B(\buu, r))}{g(r )} \geq 1 $ for any $\buu \in A$, then $  \cH_g (A)  \leq  C \mu (A) $. 
\end{itemize}
\end{proposition}
This is a standard result in Euclidian spaces: see Lemmas 2 and 3 of Rogers and Taylor \cite{RoTa} for the original proof. We refer to \cite{Edgar07} for a general version in metric spaces (more precisely, we 
refer to Theorem 4.15 \cite{Edgar07} in combination with Proposition 4.24 \cite{Edgar07}). 
We shall actually need the following more specific result.   
\begin{lemma}
\label{exactdens}
Let $\mu \in \cM_f (\boU)$. Let $g: (0, r_0)\rightarrow (0, \infty)$ be a right-continuous and non-decreasing application such that $\lim_0 g=0$ and such that there exists $C>1$ that satisfies $g(2r)\leq C g(r)$, for any $r\in (0, r_0/2)$. First assume that there is $\kappa \in (0, \infty)$ such that the following holds true. 
\begin{equation}
\label{exact}
{\rm For} \; \mu \! {\rm -almost} \; {\rm all} \; \buu \; , \quad \limsup_{r \rightarrow 0} \frac{\mu (B(\buu, r))}{g(r)} = \kappa  \; .
\end{equation}
Next, assume that there is $\kappa_0 \in (0, \kappa )  $ such that:  
\begin{equation}
\label{support}
 \cH_g \left( \left\{  \buu \in  \boU \; : \;  \limsup_{r \rightarrow 0} g(r)^{-1}\mu (B(\buu, r)) \leq \kappa_0  \right\} \right) = 0 
 \end{equation}
Then, $  \cH_g ( \, \cdot \,  \cap {\rm supp } \, \mu ) =\kappa^{-1} \mu $, 
where ${\rm supp } \, \mu$ stands for the topological support of $\mu$. 
\end{lemma}
\noi

\begin{remark}
\label{exactdensityrem}
Such a result holds true thanks to the specific properties of the balls of $(\boU, \delta )$; this can be extended to general Polish spaces if  $\mu$ satisfies a 
Strong Vitali Covering Property: see Edgar \cite{Edgar07} for a discussion of this topic. Lemma \ref{exactdens} is probably known, however the author is unable to find a reference. That is why a brief proof is provided below. \cq 
\end{remark}

\noi
{\bf Proof of Lemma \ref{exactdens}:} recall from (\ref{Budef}) notation $B_u$, $u \in \bbU$; if $A $ is a non-empty subset of $\boU$ that is not reduced to a point, then there exists $u \in \bbU$ such that $A \subset B_u $ and ${\rm diam} (B_u) = e^{-|u|}= {\rm diam} (A)$. This allows to take the set of balls as the covering set in the definition of $\cH_g$, which then coincides with the so-called spherical $g$-Hausdorff measure. Let us denote by $E$ the set of all $\buu \in \boU$ where the limsup in (\ref{exact}) holds. Then (\ref{exact}) implies that $\mu (\boU \backslash E)= 0$. Let $K $ be a compact subset of $  {\rm supp } \, \mu $. Since two balls of $\boU$ are either disjoint or one contains the other, then for any $p \geq 1$, there exists a finite sequence of pairwise disjoint balls $B_{u^p_1}$, ... , $  B_{u^p_{n_p}}$  with respective diameters $r_i^p = \exp (-|u^p_i|)  \leq p^{-1}$, $ 1 \leq i \leq n_p$, such that 
$$\forall 1 \leq i \leq n_p \; , \quad  K \cap B_{u^p_i} \neq \emptyset \; , \quad   K \subset \bigcup_{i=1}^{n_p} B_{u^p_i} \quad {\rm and} \quad \lim_{p\rightarrow \infty}  \; 
\sum_{i=1}^{n_p} g(r_i^p) = \cH_g (K) \; .$$
Since $K \subset \supp \mu$ and $K \cap B_{u^p_i} \neq \emptyset $, then $\mu (  B_{u^p_i}) >0 $, and it makes sense to define a function $f_p$ on $\boU$ by 
$$f_p(\buu)= \sum_{i=1}^{n_p} \un_{B_{u^p_i}} (\buu) \frac{g(r_i^p)}{\mu (B_{u^p_i})} \; . $$ 
Observe that $\int f_p d \mu = \sum_{i=1}^{n_p} g( r_i^p) $. Next, for any $\buu \in K \cap E$, 
denote by $j^p (\buu)$ the unique index $i \in \{ 1, \ldots , n_p \}$ such that $\buu \in B_{u^p_i}$. 
Observe that  $ f_p (\buu)$ is equal to $g (r^p_{j^p (\buu)} )/ \mu ( B ( \buu ,  r^p_{j^p (\buu)} ))$. Moreover, we have 
$$ \forall \buu \in K \cap E \; , \quad \liminf_{p \rightarrow \infty} \frac{g (r^p_{j^p (\buu)} )}{\mu ( B ( \buu ,  r^p_{j^p (\buu)} ))} \geq \liminf_{r \rightarrow 0} \frac{g(r)}{\mu (B(\buu, r))} = \kappa^{-1} \; .$$
Then, 
Fatou's lemma implies 
$$  \kappa^{-1} \mu (K \cap E) \leq \int_{K \cap E} \! \! 
\liminf_{p \rightarrow \infty} f_p (\buu)\,  \mu (d\buu) \leq 
 \liminf_{p \rightarrow \infty } \int_{\boU} \! \!  f_p (\buu) \mu ( d \buu) = \cH_g ( K) \; , $$
which entails 
\begin{equation}
\label{boundone}
\kappa^{-1} \mu (K ) \leq \cH_g ( K) \; .
\end{equation}
Let us prove the converse inequality. Thanks to the special properties of the balls of $(\boU, \delta)$, for any $\eta \in (0, 1)$, we can find 
a sequence of pairwise disjoint balls $(B( \buu_n, r_n)\, ;\, n \in \bN) $ whose diameters are smaller than  $\eta $, such that $\buu_n \in K \cap E$ for any $n \in \bN$ and such that the following holds true: 
$$ K\cap E  \, \subset \,  \bigcup_{n \in \bN} B( \buu_n, r_n) \quad {\rm and} \quad g(r_n) \leq  (1+ \eta ) \kappa^{-1} \, \mu ( B( \buu_n, r_n) ) \; , \; n \in \bN.$$
This implies that for any $\eta \in (0, 1) $, 
\begin{equation}
\label{ineqhaus}
 \cH^{(\eta)}_g (  K\cap E ) \leq   (1+ \eta ) \kappa^{-1} \, \mu ( K^{\eta})  ,  
 \end{equation}
where $K^{\eta} = \{ \buu \in \boU \; : \; \delta (\buu , K) \leq \eta \}$. 
By letting $\eta $ go to $0$, standard arguments entail 
\begin{equation}
\label{boundtwo}
 \cH_g (  K\cap E ) \leq   \kappa^{-1} \, \mu ( K)  \; .
\end{equation}
We next have to prove that $\cH_g ( {\rm supp} \, \mu  \cap (\boU \backslash E) \,  ) = 0$:  
let $b>a >0$; we set $ E_{a,b}= \{ \buu \in {\rm supp} \, \mu  \; : \;  \limsup_{r \rightarrow 0} g(r)^{-1}\mu (B(\buu, r)) \in [a, b) \}$. Proposition \ref{comparison} $(ii)$ implies that $ \cH_g ( E_{a, b} ) \leq Ca^{-1} \mu (E_{a, b})$. Therefore $\cH_g (E_{a, b})=0$ if $a > \kappa $ or if $0< a <b < \kappa$ and (\ref{support}) easily entails $\cH_g ( {\rm supp} \, \mu  \cap (\boU \backslash E) \,  ) = 0$. 
This, combined with (\ref{boundone}) and (\ref{boundtwo}), entails that for any compact set $K$,  
$\cH_g ( K \cap {\rm supp} \, \mu ) = \kappa^{-1} \mu ( K) $ and standard arguments complete the proof.   \cqfd

\section{Size biased trees.}
\label{sizebiasedsec}
We now introduce random trees with a distinguished infinite line of descent that are called size-biased trees and that have been introduced by R. Lyons, R. Pemantle and Y. Peres in \cite{LyoPemPer} to provide a simple proof of Kesten-Stigum theorem. To that end, let us first set some notation: let $T \in \bbT$ be such that $\partial T \neq \emptyset$.  Let $\buu \in \partial T$. We set 
$$ {\rm Gr (T, \buu )} = \left\{ v \in T\backslash \{ \varnothing \} 
 \; : \;  v \neq \buu_{|_{|v|}} \; {\rm and} \;  v_{|_{|v|-1}} = \buu_{|_{|v|-1}}  \; \right\} \; .$$
${\rm Gr (T, \buu )}$ is the set of individuals of $T$ whose parent belongs to the infinite line of descent determined by $\buu$ but who don't not themself belong to the $\buu$-infinite line of descent (namely, such individuals have a sibling on the $\buu$-infinite line of descent). In other words, ${\rm Gr (T, \buu )}$ is the set of the vertices where are grafted the subtrees stemming from the $\buu$-infinite line of descent.

  Let $\xi$ be an offspring distribution that satisfies Hyp($\xi$). The size-biased offspring distribution $\xi$  is the probability measure $\widehat{\xi}$ given by $\widehat{\xi} (k)= k\xi (k)/\ovvm$, $k \geq 0$. We define 
a probability distribution $\rho$ on $\bN^* \times \bN^*$ by 
$$ \rho (k , \ell ) = \un_{ \{ k \leq \ell \}} \ovvm^{-1} \xi (\ell ) = \un_{ \{ k \leq \ell \}}  \ell^{-1} \widehat{\xi} (\ell )  \; , \quad k, \ell \geq 1 \; , $$
and we call $\rho$ the repartition distribution associated with $\xi$. 
\begin{definition}
\label{sizebiasdef}
Let $(\cT^*, \bU^*): \Omega \rightarrow \bbT \times \boU$ be a $(\cF, \cG \otimes \cB (\boU) )$-measurable application; $(\cT^*, \bU^*)$ is a $\xi$-size-biased Galton-Watson tree 
(a $\widehat{ {\rm GW} } (\xi)$-tree for short) iff the following holds. 
\begin{itemize}
\item{{\rm Size-bias(1):}} for any $n \geq 0$, $\bU_{|n}^* \in \cT^*$. Moreover, if we set 
$\bU^*= ( I^*_n\, ; \, n \in \bN^*)$, then the sequence of $\bN^*\times \bN^*$-valued random variables 
$(I^*_n\, ; \,  k_{\bU_{| n-1}^* } (\cT) \, )$, $n \geq 1$ is i.i.d. with distribution $\rho$. 

\item{{\rm Size-bias(2):}} conditional on the sequence $ (\, (I^*_n\, ; \,  k_{\bU_{| n-1}^* } (\cT) \, )   \, ; \, n \geq 1)$, the subtrees $\theta_{ u } \cT $, where $u$ ranges in ${\rm Gr}(\cT^*, \bU^*) $,  are i.i.d. 
${\rm GW} (\xi )$-trees. \cq 
\end{itemize}
\end{definition}
The size-biased tree $\cT^*$ can be informally viewed as the family-tree of a population containing two kinds of individuals: the mutants and the non-mutants; each individual has an independent offspring: the non-mutants's one is distributed according to $\xi$ and the mutants's one according to $\widehat{\xi}$; moreover, 
the ancestor is a mutant and a mutant has 
exactly one child who is a mutant (its other children being non-mutants); the rank of birth of the mutant child is chosen uniformly at random among the progeny of its mutant genitor; so there is only one mutant per generation and $\bU^*$ represents the ancestral line of the mutants.  

   Size-biased trees have been introduced Lyons, Pemantle and Peres \cite{LyoPemPer}. In this paper, the authors mention related constructions: we refer to their paper for a detailed bibliographical account.
%
When $\ovvm \leq 1$, the size-biased tree $\cT^*$ has one single infinite line of descent and it is called a sintree after Aldous \cite{Alfringe}. Such a biased tree is related to (sub)critical  Galton-Watson trees conditioned on non-extinction: see Grimmett
\cite{Grimmett80}, Aldous and Pitman \cite{AlPittree98} and also \cite{Du4} for more details.

 The name "size-biased tree" is explained by the following elementary result (whose proof is left to the reader): let 
$G_1: \bbT \times \bbU \rightarrow [0, \infty) $ and $G_2: \bbT \rightarrow [0, \infty)$ be two measurables applications. Then, for any $n \geq 0$, we have  
\begin{equation}
\label{folklore}
\bE \left[ \sum_{ u \in \cT :  |u| = n } G_1 \left( \cut_u \cT \, ; \, u \right) G_2 (\theta_u \cT )  \right]
= \ovvm^n \, \bE \left[ G_1 \left( \cut_{\bU_{|n}^*} \cT^* \, ; \, \bU_{|n}^* \right)  \right]  \, \bE \left[ G_2 (\cT) \right],  
\end{equation}
with the convention that a sum other an empty set is null. This immediately entails 
$$ \bP \left( \cT_{|n}^*  \in dT \right) = \frac{Z_n (T)}{\ovvm^n} \, \bP \left( \cT_{|n} \in dT \right) \; , $$
which explains the name "size-biased" tree. Recall notation $M$ for the branching measure. We derive from (\ref{folklore}) the following key-lemma. 
\begin{lemma}
\label{keylemma}
Let $\xi$ be an offspring distribution that satisfies Hyp($\xi$). Let $(\cT^*, \bU^*)$ be a $\widehat{ {\rm GW} } (\xi)$-tree and let $\cT$ be a $ {\rm GW} (\xi)$-tree. Let $G : \bbT \times \boU \rightarrow [0, \infty)$ be $\cG \otimes \cB (\boU)$-measurable. Then, 
\begin{equation}
\label{keyformula}
\bE \left[ \int_{\boU} M (d\buu) \, G( \cT ;  \buu) \right]= \bE \left[ \, G( \cT^*;  \bU^*)\,  \right] . 
\end{equation}
\end{lemma} 
\noi
{\bf Proof:} let us first assume that $G ( \cT ; \buu)= G(  \cT_{|n} ; \buu_{|n} )$. Recall that 
$M(\theta_u \cT)= m^{-n} W_u$. Then, observe that 
$$ \int_{\boU} M (d\buu) \, G( \cT ; \buu) = \sum_{ u \in \cT : |u| = n  } G(  \cT_{|n} , \buu_{|n} ) m^{-n} W_u \; .$$
Since $\bE [W_u \; | \; u \in \cT ] = 1$,  (\ref{folklore}) implies (\ref{keyformula}). Recall notation $A_u$, $u\in \bbU$ from (\ref{defineAu}). The previous result implies that (\ref{keyformula}) holds true for applications $G$ of the form $\un_{C}$ where $C$ belong to $\cP$: 
$$ \cP = \{   (A_{u_1} \cap \ldots \cap A_{u_p} ) \times B_v \; ; \; p \in \bN^* \, , \, v, u_1, \ldots , u_p \in \bbU \} \; , $$
which is a pi-system generating $\cG \otimes \cB (\boU)$. Then, a standard monotone-class argument entails the desired result. \cqfd 

\medskip 

Let us apply Lemma \ref{keylemma} to our purpose: let $\cT$ be a ${\rm GW} (\xi)$-tree and let $\buu \in \partial \cT$;  recall notation ${\rm Gr} (\cT, \buu)$ from the beginning of the section. It is easy to prove that for any $r\in (0, 1]$, 
\begin{equation}
\label{ancestraldec}
M (B( \buu , r) ) =   \sum_{ p \geq 0}  \; \sum_{\substack{ v \in {\rm Gr} (\cT, \buu)\\ |v|=n(r)+ 1+ p }} \ovvm^{-p-n(r)-1} W_u \; , 
\end{equation}
where we recall that $n(r)= \lfloor -\log (r) \rfloor +1$. Let $(\cT^*, \bU^*)$ be a $\widehat{ {\rm GW} } (\xi)$-tree. For any $v \in {\rm Gr} ( \cT^*,\bU^*) $, we set 
$$ W^*_v = \limsup_{n\rightarrow \infty} \frac{Z_n ( \theta_v \cT^*) }{\ovvm^n} \; .$$
By Size-bias(2), observe that conditional on the $\bN^*\times \bN^*$-valued sequence 
$(\, (I^*_n\, ; \,   k_{\bU_{|n-1}^* } (\cT^*) \, )  \, ; \, n \geq 1)$, the random variables $(W^*_v \, ; 
v \in {\rm Gr} (\cT^*, \bU^*)\, ) $ are i.i.d. with the same distribution as $W$. 
Next, for any $n\geq 1$, we set 
$$ Y_n = \sum_{\substack{v \in {\rm Gr}(\cT^*, \bU^*)  \\ |v|= n }} W^*_v \quad {\rm and} \quad  
X_n = \sum_{p \geq 0} m^{-p} Y_{p+n}  \; .$$
Lemma (\ref{keylemma}) combined with (\ref{ancestraldec}) entails that 
\begin{equation}
\label{ballconsq}
\bE \left[ \int_{\boU} \! \! \! M (d\buu) \,  G  \left( \,  M ( B (\buu , r))  \, ; \, r \in (0, 1] \, \right)  \right]= 
\bE \left[  G \left( \,  m^{-n (r)-1} X_{n(r) +1} \, ; \, r \in (0, 1] \, \right) \right] \; . 
\end{equation}
Note that $(Y_n \, ; \, n \geq 1)$ is an i.i.d. sequence of random variables; thus, the $X_n$'s have the same distribution. We provide two (rough) bounds of the tail at $\infty$ of $Y_1$ and $X_1$ that are needed in the proof section: recall notation $F$ from the introduction; the first bound is a straightforward consequence of the definitions: 
\begin{equation}
\label{OKbound}
\bP ( X_1 >x) \geq \bP (Y_1 >x ) \geq C_0 \, \bP (W > x) = C_0 \,  \exp (-F(x)) \; .
\end{equation}
where $C_0:= (1-\widehat{\xi} (1))= \bP ( k_{ \{ \varnothing \}} (\cT^*) \geq 2 ) $ is a positive constant.  Next observe that 
$$ \bE \left[ \int_{\boU} \! \! \! M (d\buu) \,  G  \left( \,  M ( B (\buu , 1)) \, \right)  \right] = \bE \left[ 
 \sum_{i=1}^{k_\emptyset (\cT) } \ovvm^{-1} W_i \, G (\ovvm^{-1} W_i ) \, \right] 
  = \bE \left[ W G( \ovvm^{-1} W) \right] . $$
Therefore (\ref{ballconsq}) with $r= 1$ entails 
$$ \bE \left[ W \, G( \ovvm^{-1} W) \right]   = \bE \left[  G \left( \,  m^{-2} X_{2}\, \right) \right] \\
=  \bE \left[  G \left( \,  m^{-2} X_{1}\, \right) \right] \; . $$
Cauchy-Schwarz inequality and the previous identity both imply the following: 
 \begin{equation}
\label{roughbound}
\bP ( X_1 >\ovvm \, x)  = \bE [ W \, \un_{ \{ W> x \}}] \leq   C_1   \exp \left( -  F(x)/2 \right) \; , 
\end{equation}
where we have set  $C_1  := \sqrt{\bE \left[  W^{2}\right] }$, which is finite if $F$ satisfies (\ref{conseqF}).

\section{Proof of Theorem \ref{mainres}.}
\label{proofsec}
First observe that (\ref{OKbound}) implies that $\bP ( Y_n > F^{-1} (\log n )) \geq C_0n^{-1}$. Since the $Y_n$'s are independent, the converse of Borel-Cantelli lemma entails that $Y_n $ is larger than 
$F^{-1} (\log n )$ for infinitely many $n$ almost surely. Now, since $X_n \geq Y_n$, we easily get 
\begin{equation}
\label{lowerdens}
\bP{\rm -a.s.} \quad \limsup_{n \rightarrow \infty} \frac{m^{-n}X_n}{g(e^{-n})} \geq 1 \; ,  
 \end{equation}
where we recall that $g$ is given by (\ref{defgg}). Next (\ref{roughbound}) implies that 
$\bP ( X_n > \ovvm \, F^{-1} (3\log n )) \leq  C_1 \cdot n^{-3/2}$. Then, Borel-Cantelli combined with (\ref{recipwH}) gives 
\begin{equation}
\label{upperdens}
\bP{\rm -a.s.} \quad \limsup_{n \rightarrow \infty} \frac{m^{-n}X_n}{g(e^{-n})} \leq 6^{a} \ovvm  \; .  
\end{equation}
Kolmogorov $0$-$1$ law combined with (\ref{lowerdens}) and (\ref{upperdens})  imply that 
$\limsup_{r \rightarrow 0} g(r)^{-1}m^{-n(r)-1}X_{n(r)+1}$ is a.s. equal to a constant $\kappa =c_{\xi}^{-1} $ that is {\it positive} and {\it finite}, which entails (\ref{typicaluppdens}) by (\ref{ballconsq}).  
%

\medskip 

  Let us prove (\ref{exactH}): by Lemma \ref{exactdens}, we only need to prove that there exists $\kappa_0 >0$ such that  
\begin{equation}
\label{supportM} 
\cH_g \left( \left\{ \buu \in \partial \cT\;  : \quad  \limsup_{n \rightarrow \infty} g(r)^{-1} M (B (\buu,r)) < \kappa_0 \, \right\} \right) = 0 \; .
\end{equation}
Recall that $ M (B (\buu,r))= m^{-n(r)}W_{\buu_{|n(r)}}$. Thus, if we set 
$$E_{n_0}=  \left\{ \buu \in \partial \cT \; : \;  \forall n \geq n_0 \; , \; W_{\buu_{|n}} < \ovvm^{-1} \, F^{-1} (\frac{_1}{^2}\log n ) \, \right\}$$
then, an easy argument using (\ref{recipwH}) entails that Claim (\ref{supportM}) is true with $\kappa_0 =4^{-a} \ovvm^{-1}$ as soon as almost surely the following holds:  
\begin{equation}
\label{transl} 
\forall n_0 \geq 2 \; , \quad \cH_g \left(E_{n_0}  \right) = 0 \; . 
\end{equation}
Let us prove (\ref{transl}): to that end, we set for any $N >n_0$: 
$$ J_{n_0, N}= \left\{ v \in \cT  \; : \; |v|=N \; ;  \; \forall  n \in \{ n_0, \ldots, N \}   \; , \; W_{v_{|n} } 
< \ovvm^{-1} \, F^{-1} (\frac{_1}{^2}\log n ) \; \right\} .$$ 
Thus, $\{B_v\, ; \, v \in J_{n_0, N} \} $ is a $(e^{-N})$-cover of  $ E_{n_0} $ and 
$ \cH_g^{(e^{-N})} \left(E_{n_0} \right) \leq g(e^{-N}) \#  J_{n_0, N} $.  We now estimate $g(e^{-N}) \#  J_{n_0, N} $ thanks to (\ref{folklore}):
\begin{eqnarray*}
\bE \left[g(e^{-N}) \#  J_{n_0, N} \right] & =&g(e^{-N}) \bE \left[ \sum_{v \in \cT : |v|=N  } \un_{ \{  
\forall  n \in \{ n_0, \ldots, N \}   \; , \; W_{v_{|n} } < \ovvm^{-1} \, F^{-1} (\frac{_1}{^2}\log n )     \} } \right] \\
& =& F^{-1} (\log N) \, \bP\left(  \forall  n \in \{ n_0, \ldots, N \}   \; , \; X^*_n  <  \ovvm^{-1}\,  F^{-1} (\frac{_1}{^2}\log n )  \right) \; , 
\end{eqnarray*}
where $X^*_n = m^{-(N-n)} W' + \sum_{p=n+1}^{N} m^{-(p-n)} Y_p$ and where $W'$ is distributed as $W$ and is independent of the $Y_n$'s. Observe that $X^*_n \geq  \ovvm^{-1}Y_{n+1}$ for any $n_0 \leq n < N$. Therefore, (\ref{OKbound}) implies  
\begin{eqnarray*}
\bP \left( \forall  n \in \{ n_0, \ldots, N \}  \; , \;  X^*_n  < \ovvm^{-1} F^{-1} (\frac{_1}{^2}\log n )  \right)& \leq  &
\prod_{n=n_0}^{N-1} \bP\left( Y_{n+1} < F^{-1} (\frac{_1}{^2} \log n )  \right) \\
& \leq  &
\prod_{n=n_0}^{N-1} \left( 1-C_0 n^{-1/2}  \right)  \\
& \leq & \exp \left( -\frac{_1}{^2}C_0  \left( \sqrt{N}-\sqrt{n_0} \,  \right)  \right) \; . 
\end{eqnarray*}
Since (\ref{conseqF}) implies $F^{-1} (\log N) \leq  F^{-1} (1) (2\log N)^a$, we easily get 
 $\lim_{N\rightarrow \infty} g(e^{-N}) \#  J_{n_0, N}= 0$, which completes the proof of (\ref{transl}) and thus of (\ref{exactH}). \cqfd  
%

%


%



%
\end{document}